\documentclass[12pt]{amsart}
\usepackage{amssymb, amscd, amsmath, amsthm, latexsym, enumerate}
\newtheorem{theorem}{Theorem}

\newtheorem{proposition}[theorem]{Proposition}

\newtheorem{cor}[theorem]{Corollary}

\begin{document}

\title[Pro-$p$ groups of positive deficiency]{\boldmath Pro-$p$ groups of positive deficiency}

\author[J. A. Hillman]{Jonathan A. Hillman}
\address{School of Mathematics and Statistics\\ University of Sydney\\  Sydney, NSW 2006\\ AUSTRALIA}
\email{jonh@maths.usyd.edu.au}

\author[A. Schmidt]{Alexander Schmidt}
\address{NWF-I Mathematik\\ Universit\"at Regensburg\\ D-93040 Regensburg\\ GERMANY}
\email{alexander.schmidt@mathematik.uni-regensburg.de}

\begin{abstract}
Let $\Gamma$ be a finitely presentable pro-$p$ group with
a nontrivial finitely generated closed normal subgroup $N$ of infinite index.
Then $\mathrm{def}(\Gamma)\leq1$,
and if $\mathrm{def}(\Gamma)=1$ then $\Gamma$ is a pro-$p$
duality group of dimension~$2$,
$N$ is a free pro-$p$ group and $\Gamma/N$ is virtually free.
In particular, if the centre of $\Gamma$ is nontrivial and
$\mathrm{def}(\Gamma)\geq 1$,  then $\mathrm{def}(\Gamma)=1$, $cd\ G\leq 2$ and $\Gamma$ is
virtually a direct product $F\times\mathbb{Z}_p$,
with $F$ a finitely generated free pro-$p$ group.
\end{abstract}



\maketitle

The deficiency $\textrm{def}(\mathcal{P})$ of a finite presentation
$\mathcal{P}$ of a group $\pi$ is the difference $g-r$,
where $g$ and $r$ are the numbers of generators and relations (respectively).
Since $\textrm{def}(\mathcal{P})$ is bounded above by the
rank of the abelianization $\pi/\pi'$,
we may define the deficiency $\textrm{def}(\pi)$ of a finitely presentable group $\pi$ as the maximum
over all deficiencies of presentations for $\pi$.
A finite presentation $\mathcal{P}$ determines a
finite 2-complex $C(\mathcal{P})$ with $\pi_1(C(\mathcal{P}))\cong\pi$,
and $\chi(C(\mathcal{P}))=1-\textrm{def}(\mathcal{P})$.
Thus $\textrm{def}(\pi)\geq 1$ if and only if
there is a finite 2-complex $X$ with $\pi_1(X)\cong\pi$ and $\chi(X)\leq0$.
This property is inherited by subgroups of finite index,
since the Euler characteristic is multiplicative.
In conjunction with other group-theoretic hypotheses, this
often leads to strong constraints on the group.
For instance, if $\mathrm{def}(\pi)\geq1$
and $\beta_1^{(2)}(\pi)=0$ then $\mathrm{def}(\pi)=1$\linebreak and
$cd\ \pi=2$ or $\pi\cong\mathbb{Z}$.
(See Theorem 2.5 of \cite{Hi}.)
The $L^2$-Betti number condition holds if $\pi$ has a finitely generated
infinite normal subgroup of infinite index, or if $\pi$ has an infinite
amenable normal subgroup.
(See Chapter 7 of \cite{Lu}.)

\medskip
We are interested in finding analogous results for pro-$p$ groups.
There is at present no good $p$-adic analogue of the von Neumann algebra,
providing an invariant with the formal properties of the $L^2$-Betti numbers
in the discrete case.
Nevertheless, we shall show that if a finitely presentable pro-$p$ group
$\Gamma$ with $\mathrm{def}(\Gamma)\geq1$
has a nontrivial finitely generated closed normal subgroup
$N$ of infinite index then $\mathrm{def}(\Gamma)=1$,
$\Gamma$ is a pro-$p$ duality group of dimension~$2$,
$N$ is a free pro-$p$ group and $\Gamma/N$ is virtually free.
In particular, if a finitely presentable pro-$p$ group
has a nontrivial finite normal subgroup then it has deficiency $\leq0$.
(We do not know whether the analogue of the latter result
holds in the discrete case.)
We derive two corollaries from the main result.
Firstly, if the centre of $\Gamma$ is nontrivial and
$\mathrm{def}(\Gamma)\geq1$, then $\mathrm{def}(\Gamma)=1$,
$cd\ \Gamma\leq 2$
and  $\Gamma$ is
virtually a direct product $F\times\mathbb{Z}_p$,
with $F$ a finitely generated free pro-$p$ group.
Secondly, if $\Gamma$ has a finitely generated abelian closed normal subgroup
$A$ and $\mathrm{def}(\Gamma)\geq1$, then
$\mathrm{def}(\Gamma)=1$,
$cd\ \Gamma\leq 2$ and  $A\cong\mathbb{Z}_p$ or
$A$ has finite index in $\Gamma$ and $A\cong{\mathbb Z}_p^2$.

\medskip
If $G$ is a pro-$p$ group let
$H^i(G)=H^i(G;\mathbb{F}_p)$ and
$h^i(G)=\mathrm{dim}_{\mathbb{F}_p}H^i(G)$, where $\mathbb{F}_p$
is the trivial $G$-module of order $p$.
The following proposition is well-known.
(See, for example, Proposition 3.9.4 of \cite{NSW}.)

\begin{proposition}
Let $p$ be a prime number and let $G$ be a pro-$p$
group. Then $G$ is finitely presentable  if and only if $h^1(G) <
\infty$ and \hbox{$h^2(G) < \infty$}, and in this case
$\mathrm{def}(G)=h^1(G)-h^2(G)$.
\qed
\end{proposition}

If the numbers $h^{i}(G)$ are finite for $i=0,\ldots,n$,
we may define the $n$th partial Euler-Poincar\'{e} characteristic by
$\chi_n(G)=\sum_{i=0}^n (-1)^{i}h^{i}(G)$.
On applying Lemma 3.3.15 of \cite{NSW} (first ed.\ 3.3.12) and Shapiro's Lemma
to the (co)induced module
$A=\mathbb{F}_p[G/U]\cong \mathrm{Hom}_U(\mathbb{F}_p[G],\mathbb{F}_p)$,
we obtain the inequality
\[
\chi_n(U)(-1)^n \leq (G:U)\chi_n(G)(-1)^n
\]
for every open subgroup $U \subseteq G$.
In particular, $\chi_2(G)=1-\mathrm{def}(G)$
is submultiplicative:
if $U$ is an open subgroup of $G$ then $\chi_2(U)\leq\chi_2(G)[G:U]$.

\begin{theorem}
Let $\Gamma$ be a finitely presentable pro-$p$ group with a nontrivial
finite normal subgroup $N$. Then $\mathrm{def}(\Gamma)\leq0$.
\end{theorem}

\begin{proof}
Let $a\in N$ be an element of order $p$ and set $A=\langle a\rangle\subset N$.
Let $U=C_\Gamma(A)$ be the centralizer of $A$.
Then $U$ is a closed subgroup of finite index in $\Gamma$,
and hence is also open.
Since $A\cong\mathbb{F}_p$ is central, the group extension $0\to A \to U\to U/A\to 1$
is classified by an element $\alpha\in H^2(U/A)$.
Now $H^2(U/A)=\varinjlim{H^2}(U/V)$,
where the limit is taken over open normal
subgroups $V\leq{U}$ which contain $A$,
and so $\alpha$ is in the image of $H^2(U/V)$ for some such $V$.
Clearly the restriction of $\alpha$ to $H^2(V/A)$ vanishes.
Then $V\cong A\times{V}/A$.
As the Hochschild-Serre spectral sequence associated to a direct product
degenerates at $E^2$, see \cite{Ja} or Theorem 2.4.6 of \cite{NSW},
we obtain
$h^1(V)=h^1(V/A)+1$ and $h^2(V)=h^2(V/A) + 1 + h^1(V/A)$.
This implies $\chi_2(V)> 0$, hence $\chi_2(\Gamma)>0$
and $\textrm{def}(\Gamma)\leq 0$.
\end{proof}

In the discrete case there is no general result asserting that torsion
cohomology classes restrict to 0 on suitable finite-index subgroups.
(This is however true if the group is a surface group. See \cite{Hi02}.)

\begin{theorem}
Let $\Gamma$ be a finitely presentable pro-$p$ group with a nontrivial
finitely generated closed normal subgroup $N$ of infinite index.
Then $\mathrm{def}(\Gamma)\leq1$.
\end{theorem}

\begin{proof}
After passing to an open subgroup $U$ we may assume
that $U/N$ acts trivially on the finite group $N/N^p[N,N]$ and that the
corresponding extension splits: $U/N^p[N,N]\cong(U/N)\times{N/N^p[N,N]}$.
Hence the transgression from $H^i(U/N;H^1(N))$ to $H^{i+2}(U/N)$ in
the Hochschild-Serre spectral sequence for $U$ as an extension of $U/N$
by $N$ is trivial.
(See Theorem 2.4.4 of \cite{NSW}, first ed.\ 2.1.8.)
It follows that $h^1(U)=h^1(U/N)+h^1(N)$ and
$h^2(U)\geq{h^2(U/N)}+h^1(U/N)h^1(N)$.
Hence
\[
\begin{array}{rcl}
\chi_2(U)&\geq &1-h^1(U/N)-h^1(N)+h^1(U/N)h^1(N)\\
&=&(h^1(U/N)-1)(h^1(N)-1).
\end{array}
\]
Since $h^1(U/N)\geq1$ and $h^1(N)\geq1$ it follows that
$\chi_2(U)\geq0$ and so $\chi_2(\Gamma)\geq0$.
Therefore $\mathrm{def}(\Gamma)\leq1$.
\end{proof}

A slight sharpening of the estimate for $h^2(U)$ gives a stronger result.

\begin{theorem}
Let $\Gamma$ be a finitely presentable pro-$p$ group with
$\mathrm{def}(\Gamma)=1$ and with a nontrivial finitely generated
closed normal subgroup $N$ of infinite index.
Then $\Gamma$ is a pro-$p$ duality group of dimension~$2$,
$N$ is a free pro-$p$ group and $\Gamma/N$ is virtually free.
Moreover, either $N\cong\mathbb{Z}_p$ or $\Gamma/N$ is virtually $\mathbb{Z}_p$.
\end{theorem}

\begin{proof}
We note first that $N$ must be infinite, by Theorem 2.
If $U$ is an open subgroup of $\Gamma$ then $U\cap N$
has finite index in $N$, and thus is a nontrivial
finitely generated closed normal subgroup of infinite index in $U$.
Hence $\mathrm{def}(U)\leq1$, by Theorem 3.
Thus $0\leq\chi_2(U)\leq\chi_2(\Gamma)[\Gamma:U]=0$,
by the submultiplicativity of $\chi_2$
and the hypothesis that $\mathrm{def}(\Gamma)=1$.
Therefore $\chi_2(U)=0$ for all such subgroups $U$,
and so $cd\ \Gamma\leq2$, by Theorem 3.3.16 of \cite{NSW} (first ed.\  3.3.13).

As in Theorem 3 there is an open subgroup $U$
containing $N$ such that
$U/N^p[N,N]\cong(U/N)\times{N/N^p[N,N]}$.
Let $d_3:H^0(U/N;H^2(N))\to{H^3(U/N)}$ be the $d_3^{02}$ differential
of the Hochschild-Serre spectral sequence, and let
$c=\mathrm{dim}_{\mathbb{F}_p}\mathrm{Ker}(d_3)$.
Then $h^1(U)=h^1(U/N)+h^1(N)$ and $h^2(U)={h^2(U/N)}+h^1(U/N)h^1(N)+c$.
Since $\chi_2(U)=0$ it follows that
\[
(h^1(N)-1)(h^1(U/N)-1)+h^2(U/N)+c=0.
\]
Since $N$ and $U/N$ are each nontrivial $h^1(N)\geq1$ and $h^1(U/N)\geq1$.
Thus the three summands are all non-negative and so must be 0.
Since $h^2(U/N)=0$ the quotient $U/N$ is a free pro-$p$ group.
In particular, $\Gamma/N$ is virtually free.
Since $H^3(U/N)=0$ and $c=0$,
it follows that $H^0(U/N;H^2(N))=0$.
Since $U/N$ is a pro-$p$ group and $H^2(N)$ is a discrete $p$-torsion module,
it follows that $H^2(N)=0$.
(See Corollary 1.6.13 of \cite{NSW}, first ed.\  1.7.4.)
Thus $N$ is also a free pro-$p$ group. (The fact that $N$ is free if $h^1(U/N)=1$ follows also from \cite{Ko06},
since $\mathrm{def}(U)=1$ and $N$ is finitely generated.)

In particular, $U$ is an extension of finitely generated free pro-$p$ groups.
Hence $U$ is a pro-$p$ duality group of dimension~$2$,
by \cite{Pletch} Theorem 3.9.\footnote{The referee points out that the proofs of the main results in \cite{Pletch} are correct only for pro-$p$ groups, but not for general profinite groups.}
Since $U$ is an open subgroup of $\Gamma$ and $cd\ \Gamma\leq 2$,
the group $\Gamma$ is also a pro-$p$ duality group of dimension $2$,
by \cite{Pletch} Theorem 3.8.

If $h^1(N)=1$ then $N\cong\mathbb{Z}_p$;
otherwise $h^1(U/N)=1$ and so $\Gamma/N$ is virtually $\mathbb{Z}_p$.
\end{proof}

\begin{cor} \label{cor1}
Let $\Gamma$ be a finitely presentable pro-$p$ group with nontrivial centre $\zeta\Gamma$.
Then $\mathrm{def}(\Gamma)\leq1$.
If\/ $\mathrm{def}(\Gamma)=1$ then either $\Gamma\cong{\mathbb{Z}}_p^2$
or $\zeta\Gamma\cong{\mathbb{Z}}_p$.
Moreover, $\Gamma$ is virtually a direct product $F\times\mathbb{Z}_p$,
with $F$ a finitely generated free pro-$p$ group.
\end{cor}

\begin{proof}
We may clearly assume that $\mathrm{def}(\Gamma)\geq1$,
since there is nothing to prove if $\mathrm{def}(\Gamma)<1$.

If $[\Gamma:\zeta\Gamma]$ is finite the commutator subgroup
$\Gamma'=[\Gamma,\Gamma]$ is finite, by a lemma of Schur.
(See Proposition 10.1.4 of \cite{Rob}.
The argument given there for the discrete case
extends without change to the pro-$p$ case.)
Hence $\Gamma'=1$, by Theorem 2, and so $\Gamma$ is abelian.
Moreover, $\Gamma$ is torsion-free, by Theorem 2 again.
Therefore $\Gamma\cong \mathbb{Z}_p^2$ or $\mathbb{Z}_p$,
since $\mathrm{def}(\Gamma)\geq1$.
In particular,
$\Gamma\cong{F}\times\mathbb{Z}_p$ with $F$ free of rank 1 or 0.

Suppose now that $[\Gamma:\zeta\Gamma]=\infty$.
Let $C$ be a nontrivial finitely generated closed subgroup of $\zeta\Gamma$.
Then $C$ is a closed normal subgroup of infinite index in $\Gamma$,
and is infinite, by Theorem 2.
Hence $\mathrm{def}(\Gamma)=1$,
$C$ is free and $\Gamma/C$ is virtually free, by Theorem 3.
But then $C\cong\mathbb{Z}_p$ and $\zeta\Gamma/C$ is finite,
since it is a central closed subgroup of infinite index in a virtually
free pro-$p$ group.
Thus $\zeta\Gamma$ is also finitely generated and so
$\zeta\Gamma\cong\mathbb{Z}_p$, by the same argument.

Let $F$ be a free pro-$p$ subgroup of finite index in $\Gamma/\zeta\Gamma$.
Since $\Gamma/\zeta\Gamma$ is virtually free $\Gamma$ has an open subgroup $U$
containing $\zeta\Gamma$ and such that
$F=U/\zeta\Gamma$ is a finitely generated free pro-$p$ group.
Since $F$ is free and acts trivially on $\zeta\Gamma$ the extension splits,
and so $U\cong{F}\times\mathbb{Z}_p$.
\end{proof}

Such direct products $F\times\mathbb{Z}_p$ with
$F$ a finitely generated free pro-$p$ group of rank $>1$
clearly have deficiency 1 and centre $\mathbb{Z}_p$.

\begin{cor}
Let $\Gamma$ be a finitely presentable pro-$p$ group with a nontrivial
finitely generated closed abelian normal subgroup $A$.
Then we have $\mathrm{def}(\Gamma)\leq1$.
If\/ $\mathrm{def}(\Gamma)=1$, then $cd\ \Gamma \leq 2$. Moreover,
$A\cong {\mathbb Z}_p$ or $A$ has finite index and $A\cong{\mathbb Z}_p^2$.
\end{cor}

\begin{proof}
If $A$ has a nontrivial torsion subgroup $T$,
then $T$ is finite,
and so $\mathrm{def}(\Gamma)< 1$, by Theorem~2.
So we may assume that $A \cong{\mathbb Z}_p^n$ for some $n\geq 1$.
If $A$ has infinite index in $\Gamma$,
then $\mathrm{def}(\Gamma)\leq 1$ by Theorem~3.
So assume that $(\Gamma:A)<\infty$.
Then $0\leq\chi_2(A)\leq (\Gamma:A)\chi_2(\Gamma)$, implying $\chi_2(\Gamma)\geq 0$, hence $\mathrm{def}(\Gamma)\leq 1$. The same argument shows $\mathrm{def}(\Gamma)< 1$ if $n\geq 3$.

Finally, assume that $\mathrm{def}(\Gamma)=1$. If $A$ has finite index in $\Gamma$, then $A\cong {\mathbb Z}_p$ or $A\cong {\mathbb Z}_p^2$  by the argument above. The open subgroups  $U$ of $A$ are cofinal among the open subgroups of $\Gamma$ and we have $\chi_2(U)=0$. Since $\chi_2(\Gamma)=1-\mathrm{def}(\Gamma)=0$, we obtain $cd\ \Gamma\leq 2$ by Theorem 3.3.16 of \cite{NSW} (first ed.\  3.3.13).

If $A$ has infinite index, it is free by Theorem~4. Hence $A\cong {\mathbb Z}_p$ in this case. Furthermore $cd\ \Gamma=2$, again by Theorem~4.
\end{proof}

\medskip
The questions considered here can be traced back to Murasugi and Gottlieb.
Murasugi \cite{Mu65} conjectured that if a finitely presentable (discrete)
group $\pi$ has nontrivial centre then $\mathrm{def}(\pi)\leq1$,
with equality only if $\pi\cong {\mathbb Z}^2$ or $\zeta\pi\cong {\mathbb Z}$,
and verified this for one-relator groups and classical link groups,
while Gottlieb \cite{Go65} showed that if the fundamental group
of an aspherical finite complex $X$ has nontrivial centre then $\chi(X)=0$.

\medskip
Nakamura \cite{Na} has given a direct analogue of Gottlieb's Theorem
for pro-$p$ groups:
if $G$ is a pro-$p$ group with $\zeta{G}\not=1$,
$cd_pG<\infty$ and $\beta_i(G;\mathbb{F}_p)<\infty$
for all $i$ then $\chi(G;\mathbb{F}_p)=0$.
(He uses a pro-$p$ version of Stallings' argument involving
universal trace functions.)
This result and the above Corollary~\ref{cor1} overlap, but neither implies the other.

\medskip
Amenable groups have no nonabelian free subgroups.
The latter notion extends naturally to the pro-$p$ case.
If a finitely presentable discrete group has deficiency $>1$ then it
contains a nonabelian free group \cite{Rom}.
(In fact every such group is ``large": it has a subgroup of finite index which
maps onto a free group of rank 2 \cite{BP}.)
If a finitely presentable pro-$p$ group $\Gamma$ with $h^1(\Gamma)=d>1$
has no nonabelian free pro-$p$ subgroup then $h^2(\Gamma)\geq d^2/4$ by \cite{Zel}.
Thus either $d=2$ and $\Gamma$ is a one-relator pro-$p$ group
or $\mathrm{def}(\Gamma)\leq0$.
Do either of these arguments extend to (discrete or pro-$p$)
groups with infinite normal subgroups having no nonabelian free subgroup?

\vskip1cm


\begin{thebibliography}{NSW2}

\bibitem[BP]{BP} Baumslag, B. and Pride, S.J. Groups with two more generators than
relators,
J. London Math. Soc. 17 (1978), 425--426.

\bibitem[Go]{Go65} Gottlieb, D. H. A certain subgroup of the fundamental group,
Amer. J. Math. 87 (1965), 840--856.

\bibitem[Hi1]{Hi} Hillman, J. A. {\em Four-Manifolds, Geometries and Knots}, GT Monograph vol. 5,
Geometry and Topology Publications, Coventry, University of Warwick 2002. Revision 2007.

\bibitem[Hi2]{Hi02} Hillman, J. A. Deficiencies of lattices in connected Lie groups, Bull. Austral. Math. Soc. 65 (2002), 393--397.

\bibitem[Ja]{Ja} Jannsen, U. The splitting of the Hochschild-Serre spectral sequence for a pro\-duct of groups, Canad.
Math. Bull. 33 (1990) 181--183.


\bibitem[Ko]{Ko06} Kochloukova, D. H. On a conjecture of E. Rapaport Strasser about knot-like groups and its pro-$p$ version, J. Pure App. Alg. 204 (2006), 536--554.

\bibitem[Lue]{Lu} L\"uck, W. {\em $L^2$-Invariants: Theory and Applications
to Geometry and $K$-Theory},
Ergebnisse 3. Folge, Bd. 44,
Springer-Verlag Berlin, Heidelberg, New York 2002.

\bibitem[Mu]{Mu65} Murasugi, K. On the centre of the group of a link,
Proc. Amer. Math. Soc. 16 (1965), 1052--1057.

\bibitem[Na]{Na} Nakamura, H. On the pro-$p$ Gottlieb Theorem,
Proc. Japan Acad. Math. Sci. 68 (1992), 279--282.

\bibitem[NSW$^2$]{NSW} Neukirch, J., Schmidt, A. and Wingberg, K. {\em Cohomology of
Number Fields, 2nd ed.}, Grundlehren Bd. 323, Springer-Verlag Berlin, Heidelberg, New York 2008.

\bibitem[Ple]{Pletch} Pletch, A. Profinite duality groups. I,  J. Pure Appl. Algebra 16 (1980), 55--74.

\bibitem[Rob]{Rob} Robinson, D. S. {\em A Course in the Theory of Groups}, Graduate Texts in Mathematics 80, Springer-Verlag Berlin, Heidelberg, New York 1982.

\bibitem[Rom]{Rom} Romanovskii, N. S. Free subgroups of finitely-presented groups, Algebra and Logic 16 (1977), 62--68.

\bibitem[Ze]{Zel} Zelmanov, E. On groups satisfying the Golod-Shafarevich condition, in {\em New Horizons in pro-$p$ Groups} (edited by Du Sautoy, M., Segal, D., and Shalev, A.), Birkh\"auser Verlag Boston, Basel, Berlin 2000, 223--232.

\end{thebibliography}
\end{document}